\documentclass[envcountsect]{svmult}
\pdfoutput=1
\usepackage{graphicx}
\usepackage{amssymb, amsmath}
\usepackage[bookmarks=false]{hyperref}
\usepackage{caption}
\usepackage{mathptmx, helvet, courier, makeidx, graphicx, multicol}
\usepackage{mathrsfs}
\usepackage[bottom]{footmisc}

\usepackage{cite}
\newcommand\numberthis{\addtocounter{equation}{1}\tag{\theequation}}

\allowdisplaybreaks[1]

\spnewtheorem{Theorem}[theorem]{Theorem}{\bfseries}{\itshape}
\spnewtheorem{Lemma}[theorem]{Lemma}{\bfseries}{\itshape}
\spnewtheorem{Corollary}[theorem]{Corollary}{\bfseries}{\itshape}
\spnewtheorem{Proposition}[theorem]{Proposition}{\bfseries}{\itshape}
\spnewtheorem{Remark}[theorem]{Remark}{\itshape}{\rmfamily}
\spnewtheorem{Definition}[theorem]{Definition}{\bfseries}{\itshape}

\smartqed

\begin{document}
\title*{Mock Theta Function Identities Deriving from  Bilateral Basic
Hypergeometric Series}

\author{James Mc Laughlin\thanks{This work was partially supported by a grant from the Simons Foundation (\#209175 to James Mc Laughlin).}}

\institute{James Mc Laughlin \at Mathematics Department,
 25 University Avenue,
West Chester University, West Chester, PA 19383. \email{jmclaughlin2@wcupa.edu}}

\motto{This paper is dedicated to Krishna Alladi on the occasion of his $60^{th}$ birthday.}
\maketitle


\abstract{The bilateral series corresponding to many of the third-, fifth-, sixth- and eighth order mock theta functions may be derived as special cases of $_2\psi_2$ series
\[
\sum_{n=-\infty}^{\infty}\frac{(a,c;q)_n}{(b,d;q)_n}z^n.
\]
Three transformation formulae for this series due to Bailey are used to derive various transformation and summation formulae for both these mock theta functions and the corresponding bilateral series.
\\
New and existing summation formulae for these bilateral series are also used to make explicit in a number of cases the fact that for a mock theta function, say $\chi(q)$, and a root of unity in a certain class, say $\zeta$,  that there is a theta function $\theta_{\chi}(q)$ such that
\[
\lim_{q \to \zeta}(\chi(q) - \theta_{\chi}(q))
\]
exists, as $q \to \zeta$ from within the unit circle.
\\[12pt]
{\bf 2010 Mathematics Subject Classification.} {ADD MSC NUMBERS}
}

\renewcommand{\theequation}{\arabic{section}.\arabic{equation}}

\section{Introduction}
\setcounter{equation}{0}
The mock theta functions were introduced to the world by Ramanujan in his last letter to G.H. Hardy (\cite[pp. 354--355]{R00}, \cite[pp. 220--223]{BR95}), in which he also gave examples of mock theta functions of orders
three, five and seven. Ramanujan did not explain precisely what he meant by a mock theta function, and  Ramanujan's statements were interpreted by Andrews and Hickerson \cite{AH91}  to mean  a  function $f(q)$ defined by a $q$-series which converges for $|q| < 1$ and
which satisfies the following two conditions:\\
(0) For every root of unity $\zeta$, there is a $\theta$-function $\theta_{\zeta}(q)$ such that the
difference $f(q) - \theta_{\zeta}(q)$ is bounded as $q \to \zeta$ radially.\\
(1) There is no single $\theta$-function which works for all $\zeta$; i.e., for every
$\theta$-function $\theta (q)$ there is some root of unity $\zeta$ for which $f(q)-\theta (q)$ is
unbounded as $q \to \zeta$ radially.

A similar definition was given by Gordon and McIntosh \cite{GM00}, \cite{GM12}, where they also distinguish between a mock theta function and a ``strong'' mock theta function. The modern view of mock theta functions is based on the work of Zwegers \cite{Z01,Z02}, who showed that the mock theta functions are  holomorphic parts of certain
harmonic weak Maass forms.

In relation to the results in the present paper, we recall two areas of investigation in the subject of mock theta functions. Firstly, as regards condition (0) above, Folsom, Ono and Rhoades \cite{FOR13} make this condition explicit for the third order mock theta function $f(q)$, in that they found a formula for the $\theta$-function $\theta_{\zeta}(q)$ and an expression for the limit of the
difference $f(q) - \theta_{\zeta}(q)$  as $q \to \zeta$ radially, where $\zeta$ is a primitive even-order root of unity (see Theorem \ref{fort} below). Secondly, there is the subject of basic hypergeometric transformations of mock theta functions, and summation formula for sums/differences of mock theta functions. Several identities of these types were stated by Ramanujan \cite{R00}
 and were subsequently investigated by Watson \cite{W36},  and later work was carried out by  Andrews \cite{A66a}, \cite{A66c}, and more recently by Gordon and McIntosh \cite{GM03}, \cite{GM12}.

The starting point for the investigation in the present paper is the observation that many of the mock theta functions are special cases of one ``side'' ($n\geq 0$ or $n<0$) of certain general bilateral series, bilateral series which in turn derive  from the $_2\psi_2$ series
\begin{equation*}
\sum_{n=-\infty}^{\infty}\frac{(a,c;q)_n}{(b,d;q)_n}z^n
=
\sum_{n=0}^{\infty}\frac{(a,c;q)_n}{(b,d;q)_n}z^n+\sum_{n=1}^{\infty}\frac{(q/b,q/d;q)_n}{(q/a,q/c;q)_n}\left(\frac{bd}{acz}\right)^n.
\end{equation*}
These include:
\begin{align*}
&\text{Third order - all 9 (Ramanujan, Watson, Gordon and McIntosh)};\\
&\text{Fifth order - 8 of 10 (Ramanujan)};\\
&\text{Sixth order - 8 (Ramanujan)};\\
&\text{Eighth order -  4 of 8 (Gordon and McIntosh)}.
\end{align*}

A number of transformations and summation formula for the $_2\psi_2$ series due to Bailey \cite{B50} are combined with the representation of  these  mock theta functions in terms of the $_2\psi_2$ series,  together with other existing summation and transformation formulae for $q$-series, to derive new representations for the mock theta functions,  and other $q$-series  identities.

Results in the present paper include: \\
1)  radial limit results for a number of third-, fifth-, sixth- and eighth order mock theta functions similar to that of Folsom, Ono and Rhoades \cite{FOR13} alluded to above,  \\
2) new summation formulae for  the bilateral series associated with some of these order mock theta functions, \\
3) new transformation formulae for some of these mock theta functions deriving from these general bilateral transformations,  \\
4) a number of other summation formulae.

One example of a new summation formula is the following identity for the third order mock theta function $\phi(q)$:
\[
\phi(q)+\sum_{r=1}^{\infty}(-1;q^2)_rq^r=\sum_{n=-\infty}^{\infty} \frac{q^{n^2}}{(-q^2;q^2)_n}=\frac{(-q,-q,q^2;q^2)_{\infty}}
{(q,-q^2;q^2)_{\infty}}.
\]

This formula in turn implies that  if $\zeta$ is a primitive even-order $4k$ root of unity, then as $q$
approaches $\zeta$ radially within the unit disk,
\begin{equation*}
\lim_{q\to \zeta}
\left(\phi(q)- \frac{(q^2,-q,-q;q^2)_{\infty}}{(-q^2,q;q^2)_{\infty}}\right)= -2\sum_{n=0}^{k-1} (1+\zeta^2)(1+\zeta^4)\dots (1+\zeta^{2n})\zeta^{n+1}.
\end{equation*}

For the third order mock theta function $\psi(q)$ there is the transformation
\begin{equation*}
\psi(q)
=
-\sum_{n=0}^{\infty}(q;q^2)_n (-1)^{n}
+\frac{1}{2(q^2;q^2)^2_{\infty}}
\sum_{r=-\infty}^{\infty}
q^{2r^2+r}(4r+1)(-1)^r.
\end{equation*}
\emph{Note: See the remark at the end of the proof of Theorem \ref{mt3bsum1} about the convergence of the first series on the right.}

As an example of one of the new summation formulae there is the following:
{\allowdisplaybreaks
\begin{multline*}
\sum_{r=-\infty}^{\infty}(10r+1) q^{(5r^2+r)/2}=
\\
\phantom{asdasda}\left(
\frac{4q(q^4,q^{16},q^{20};q^{20})_{\infty}}{(q^2;q^{4})_{\infty}}
+
\frac{(q^2,q^{3},q^{5};q^{5})_{\infty}}{(-q;q)_{\infty}}
\right)
\frac{(q;q)^2_{\infty}}{(-q;q)_{\infty}}.
\end{multline*}}

A  number of results of a similar nature may be found throughout the paper.

Remark:  The first version of the present paper was written in 2014, and subsequently  the author attention was directed (my thanks to the anonymous referee) to a number of recent papers containing similar results, of which the present author was previously unaware.

In \cite{M14}, Mortenson derived several identities involving the Appell-Lerch sum
\begin{equation}\label{mxqzeq}
m(x,q,z):=\frac{1}{j(z;q)}\sum_{r=-\infty}\frac{(-1)^rq^{r(r-1)/2}z^r}{1-q^{r-1}xz},
\end{equation}
(here $j(z;q)=(z,q/z,q;q)_{\infty}$) and the universal mock theta function $g(x,q)$ (see \eqref{g3xq}), and some of these were used in \cite{M15} to derive explicit radial limits for mock theta functions. As well deriving such radial limits for several particular mock theta functions, in \cite{M15} the author also derives a general result for $g(x,q)$, a result which permits an explicit radial limit to be derived for any even-order mock theta function that may be expressed in terms of $g(x,q)$. We will compare results in the present paper with those in \cite{M14, M15} in several places throughout the paper.

For example, by applying a formula of Mortenson (\cite[Eq. (6.10)]{M15}), a different radial limit result is obtained for the eighth order mock theta function $S_0(q)$ (see \eqref{mock8radeq3}).

Subsequent to writing the first draft of the present paper, the author was also directed to the recent paper \cite{BKLMR15}, in which the authors also derive explicit radial limits for all of Ramanujan's third and fifth order mock theta functions, as well as giving the level and weight information for the theta functions (which are modular forms). The authors in \cite{BKLMR15} also state, without proof, explicit radial limits for many of the even-order mock theta functions.

In the present paper we also derive these explicit radial limits using somewhat different methods, but in addition also derive many  identities that come from the aforementioned connections with the $_2\psi_2$ series.

\section{Some required basic hypergeometric formulae}
\setcounter{equation}{0}

To prove some of the results in the present paper, it is necessary to use a number of transformation- and summation formulae for basic hypergeometric series.
{\allowdisplaybreaks
\begin{equation}\label{1def4}
(a;q)_{-n}:=\frac{(a;q)_{\infty}}{(aq^{-n};q)_{\infty}}
=\frac{1}{(aq^{-n};q)_{n}}=\frac{(-q/a)^n}{(q/a;q)_n}q^{n(n-1)/2},
\end{equation}
\begin{equation}\label{6t1eq}
\sum_{n=-\infty}^{\infty}(-z)^n q^{n^2}=(zq,q/z, q^2;q^2)_{\infty}.
\end{equation}
\begin{align}\label{q7.11eq1}
\sum_{n=-\infty}^{\infty}&\frac{(a,c;q)_n}{(b,d;q)_n}z^n
=
\frac{(az,cz,qb/acz,qd/acz;q)_{\infty}}
{(b,d,q/a,q/c;q)_{\infty}}
\\
&\phantom{asdsasdsasdsasdsasds}\times\sum_{n=-\infty}^{\infty}\frac{(acz/b,acz/d;q)_n}{(az,cz;q)_n}\left( \frac{bd}{acz} \right)^n.\notag\\
\sum_{n=-\infty}^{\infty}&\frac{(a,c;q)_n}{(b,d;q)_n}z^n\numberthis \label{q7.11eq2}\\
&=
\frac{(b/a,d/c,az,qb/acz;q)_{\infty}}
{(b,q/c,z,bd/caz;q)_{\infty}}
\sum_{n=-\infty}^{\infty}
\frac{(a,acz/b;q)_n}{(az,d;q)_n}\left( \frac{b}{a} \right)^n.\notag
\end{align}
\begin{multline}\label{8b2psi2}
\sum_{n=-\infty}^{\infty}\frac{(e,f;q)_n}
{(aq/c,aq/d;q)_n}\left(\frac{q a}{ef}\right)^n
=\frac{ (q/c,q/d,aq/e,aq/f;q)_{\infty} } {
(aq,q/a,aq/cd,aq/ef;q)_{\infty} }\\
\times \sum_{n=-\infty}^{\infty}\frac{(1-aq^{2n})(c,d,e,f;q)_n}
{(1-a)(aq/c,aq/d,aq/e,aq/f;q)_n}\left(\frac{q a^{3}}{cdef}\right)^n q^{n^2}.
\end{multline}
\begin{multline}\label{8baileyeq1}
\sum_{n=-\infty}^{\infty}\frac{(q\sqrt{a},-q\sqrt{a},b,c,d,e;q)_n}
{(\sqrt{a},-\sqrt{a},aq/b,aq/c,aq/d,aq/e;q)_n}\left(\frac{q a^2}{bcde}\right)^n\\
=\frac{ (aq,aq/bc,aq/bd,aq/be,aq/cd,aq/ce,aq/de,q,q/a;q)_{\infty} } {
(aq/b,aq/c,aq/d,aq/e,q/b,q/c,q/d,q/e,qa^2/bcde;q)_{\infty} }.
\end{multline}
\begin{multline}\label{8baileyeqc12}
\sum_{n=-\infty}^{\infty}\frac{(b,c;q)_n}
{(aq/b,aq/c;q)_n}\left(\frac{-q a}{bc}\right)^n\\
=\frac{ (aq/bc;q)_{\infty} (aq^2/b^2,aq^2/c^2,q^2,aq,q/a;q^2)_{\infty}} {
(aq/b,aq/c,q/b,q/c,-qa/bc;q)_{\infty} }.
\end{multline}
}

The identity at \eqref{6t1eq} is the famous Jacobi triple product identity.
The bilateral transformations at \eqref{q7.11eq1}, \eqref{q7.11eq2} and \eqref{8b2psi2} are all due to Bailey \cite{B50}. The identity at \eqref{8baileyeq1}  Bailey's $_6\psi_6$ summation formula  and \eqref{8baileyeqc12} is a special case of this (see \cite[Eq. (II.30), p. 357]{GR04}).

\section{  Mock theta functions of the third order}
\setcounter{equation}{0}
The third order mock theta functions stated by Ramanujan (\cite[pp. 354--355]{R00}, \cite[pp. 220--223]{BR95}) are the following basic hypergeometric series:
{\allowdisplaybreaks
\begin{align*}
f(q)&=\sum_{n=0}^{\infty} \frac{q^{n^2}}{(-q,-q;q)_n},&
\phi(q)&=\sum_{n=0}^{\infty} \frac{q^{n^2}}{(-q^2;q^2)_n},\\
\chi(q)&=\sum_{n=0}^{\infty} \frac{q^{n^2}(-q;q)_n}{(-q^3;q^3)_n},&
\psi(q)&=\sum_{n=1}^{\infty} \frac{q^{n^2}}{(q;q^2)_n}.
\end{align*}
}

All of the third order mock theta functions of Ramanujan, as well as  those stated later by Watson \cite{W36} and  Gordon and McIntosh \cite{GM03}, may be expressed in terms of the function $g(x,q)$, where
\begin{equation}\label{g3xq}
g(x,q):= \sum_{n=0}^{\infty}\frac{ q^{n^2+n}}{(x,q/x;q)_{n+1}}=
x^{-1} \left(
-1+\sum_{n=0}^{\infty} \frac{q^{n^2}}{(x;q)_{n+1}(q/x;q)_{n}}
\right).
\end{equation}
This was shown by Hickerson and Mortenson \cite[Eqs. (5.4) - (5.10)]{HM14} (this function was also defined by Gordon and McIntosh \cite{GM12}, where it was labelled ``$g_3(x,q)$''). For completeness, we consider a generalization, namely the series
\begin{equation}\label{gstqeq}
G_3(s,t,q):=1+\sum_{n=1}^{\infty}\frac{s^n t^n q^{n^2}}{(sq,tq;q)_n},
\end{equation}
which was defined in \cite[Eq. (7)]{Ch11}, and state a number of transformation formulae for this function. Note that the connection with the third order mock theta functions is that
\begin{equation}\label{13gGeq}
G_3(x,q/x,q) = (1-x)(1-q/x)g(x,q).
\end{equation}

\begin{Proposition}
Let $G_3(s,t,q)$ be as defined at \eqref{gstqeq} above. Then
{\allowdisplaybreaks
\begin{align}
G_3(s&,t,q)=-\sum_{r=1}^{\infty}
(s^{-1},t^{-1};q)_r  q^{r}+\frac{(q/s,q/t;q)_{\infty}}{(sq,tq;q)_{\infty}}
\sum_{r=-\infty}^{\infty}
(s,t;q)_r  q^{r};\numberthis \label{28t1eq4}\\
&=-\sum_{r=1}^{\infty}
(s^{-1},t^{-1};q)_r  q^{r}+\frac{(q/t;q)_{\infty}}{(sq,q;q)_{\infty}}
\sum_{r=-\infty}^{\infty}
\frac{(t;q)_r (-s)^r q^{r(r+1)/2}}{(tq;q)_r};\numberthis \label{28t1eq5}\\
&=-\sum_{r=1}^{\infty}
(s^{-1},t^{-1};q)_r  q^{r} \numberthis \label{28t1eq6}\\
&\phantom{asds}+\frac{(q/s,q/t;q)_{\infty}}{(stq,q/(st),q;q)_{\infty}}
\sum_{r=-\infty}^{\infty}
\frac{(1-stq^{2r})(s,t;q)_r (st)^{2r} q^{2r^2}}{(1-st)(sq,tq;q)_r}. \notag
\end{align}
}
\end{Proposition}
\begin{proof}
The transformations at \eqref{28t1eq4} and \eqref{28t1eq5} will follow as special cases of two more general identities. Replace $z$ with $zq/ac$, let $a, c \to \infty$ and set $b=sq$ and $d=tq$ in, respectively,  \eqref{q7.11eq1}
and  \eqref{q7.11eq2}, to get that
{\allowdisplaybreaks
\begin{align}\label{gzeq1}
\sum_{n=-\infty}^{\infty}\frac{z^n q^{n^2}}{(sq,tq;q)_n}&= \frac{(sq/z,tq/z;q)_{\infty}}{(sq,tq;q)_{\infty}}
\sum_{r=-\infty}^{\infty}
(z/s,z/t;q)_r  \left(\frac{stq}{z}\right)^{r},\\
&=\frac{(qs/z;q)_{\infty}}{(sq,stq/z;q)_{\infty}}
\sum_{r=-\infty}^{\infty}
\frac{(z/s;q)_r (-s)^r q^{r(r+1)/2}}{(tq;q)_r}.\numberthis \label{gzeq2}
\end{align}}
Lastly, replace $z$ with $st$, and use \eqref{1def4} on the terms of negative index in the new series on the left sides.

We also prove a generalization of the transformation at \eqref{28t1eq6} first, by  letting $e, f \to \infty$ in \eqref{8b2psi2}, and then replacing $a$ with $z$, $c$ with $z/s$ and $d$ with $z/t$, to get
\begin{multline}\label{gzeq3}
\sum_{n=-\infty}^{\infty}\frac{z^n q^{n^2}}{(sq,tq;q)_n}\\
=\frac{(sq/z,tq/z;q)_{\infty}}{(zq,q/z,stq/z;q)_{\infty}}
\sum_{r=-\infty}^{\infty}
\frac{(1-zq^{2r})(z/s,z/t;q)_r (zst)^{r} q^{2r^2}}{(1-z)(sq,tq;q)_r}.
\end{multline}
The identity at \eqref{28t1eq6} follows after replacing $z$ with $st$.
\qed
\end{proof}
The identities \eqref{28t1eq4} - \eqref{28t1eq6} may be more concisely expressed using the function
\begin{equation}\label{28g*eq}
G_3^{*}(s,t,q):=\sum_{n=-\infty}^{\infty}\frac{s^n t^n q^{n^2}}{(sq,tq;q)_n}
\end{equation}
as follows
\begin{align}\label{28g*eq1}
G_3^{*}(s,t,q)&=\frac{(q/s,q/t;q)_{\infty}}{(sq,tq;q)_{\infty}}G_3^{*}(s^{-1},t^{-1},q),\\
&=\frac{(q/t;q)_{\infty}}{(sq,q;q)_{\infty}}
\sum_{r=-\infty}^{\infty}
\frac{(t;q)_r (-s)^r q^{r(r+1)/2}}{(tq;q)_r};\numberthis \label{28g*eq2}\\
&=\frac{(q/s,q/t;q)_{\infty}}{(stq,q/(st),q;q)_{\infty}}
\sum_{r=-\infty}^{\infty}
\frac{(1-stq^{2r})(s,t;q)_r (st)^{2r} q^{2r^2}}{(1-st)(sq,tq;q)_r}.\numberthis \label{28g*eq3}
\end{align}
The identity at \eqref{28t1eq4} (or \eqref{28g*eq2}) was also proved by Choi \cite[Theorem 4]{Ch11}, and stated previously by Ramanujan  (see \cite[Entry 3.4.7]{ABRLNII}).

We will employ \eqref{28t1eq5} to derive some results on explicit radial limits, as mentioned earlier. Before coming to that, we remark that other transformations listed above may be used to derive some new transformations for three of the third order mock theta functions of Ramanujan and one of the third order mock theta functions of Watson (similar results may be derived for the other third order mock theta functions of Watson \cite{W36} and those of Gordon and McIntosh \cite{GM03}).
Before stating the next theorem, we recall  Watson's \cite{W36} third order mock theta function $\nu(q)$, where
\[
\nu(q)=\sum_{r=0}^{\infty}\frac{q^{n^2+n}}{(-q;q^2)_{n+1}}.
\]
\begin{Theorem}\label{mt3bsum1}
If $|q|<1$, then
\begin{equation}\label{fq3id1}
f(q)
=
-\sum_{n=1}^{\infty}(-1,-1;q)_n q^{n}+4\frac{(-q;q)^2_{\infty}}{(q;q)^3_{\infty}}
\sum_{r=-\infty}^{\infty}
\frac{ q^{2r^2+r}(4rq^r+1)}{(1+q^r)^2}.
\end{equation}
\begin{equation}\label{phiq3id1}
\phi(q)
=
-\sum_{n=1}^{\infty}(-1;q^2)_n q^{n}+4\frac{(-q^2;q^2)_{\infty}}{(q;q)^3_{\infty}}
\sum_{r=-\infty}^{\infty}
\frac{ q^{2r^2+2r}(2rq^{2r}+1)}{(1+q^{2r})^2}.
\end{equation}
\begin{multline}\label{nuq3id1}
\nu(q)
=
-\sum_{n=0}^{\infty}(-q;q^2)_n q^{n}\\
+4\frac{(-q;q^2)_{\infty}}{(q;q)^3_{\infty}}
\sum_{r=-\infty}^{\infty}
\frac{ q^{2r^2+2r}(r+1)}{(1+q^{2r+1})^2}
-2\frac{(-q;q^2)_{\infty}}{(q;q)^3_{\infty}}(-q^4,-q^{12},q^{16};q^{16})_{\infty}.
\end{multline}
\begin{equation}\label{psiq3id1}
\psi(q)
=
-\sum_{n=0}^{\infty}(q;q^2)_n (-1)^{n}
+\frac{1}{2(q^2;q^2)^2_{\infty}}
\sum_{r=-\infty}^{\infty}
q^{2r^2+r}(4r+1)(-1)^r.
\end{equation}
\end{Theorem}

\begin{proof}
For \eqref{fq3id1}, replace $z$ with $z^2$, $s$ and $t$ with $-z$ in \eqref{gzeq3}, and then let $z\to 1$.

A similar application of \eqref{gzeq3}, again with  $z$ replaced with $z^2$, $s$ replaced with $iz$ and $t$ replaced with $-iz$ and once again letting $z\to 1$ leads to \eqref{phiq3id1}.

For \eqref{nuq3id1}, in \eqref{gzeq3} again replace  $z$  with $z^2$, and then replace$s$  with $iz$, $t$  with $-iz$, let $z\to \sqrt{q}$ and divide through by $1+q$.

Finally, the transformation at \eqref{psiq3id1} follows similarly from \eqref{gzeq3}, this time with  $z$ replaced with $z^2$, $s$ replaced with $z/\sqrt{q}$ and $t$ replaced with $-z/\sqrt{q}$ and once again letting $z\to 1$. Note that convergence of the first series on the right of \eqref{psiq3id1} is in the Ces\`aro sense.
\qed
\end{proof}

As Watson pointed out in \cite[Section 7]{W37}, certain bilateral series related to fifth order mock theta functions, which are essentially the sums of pairs of fifth order mock theta functions, are expressible as theta functions, or combinations of infinite $q$-products. It seems less well known that the bilateral series associated with two of Ramanujan's third order mock theta functions are also expressible as infinite products. We also give similar statement for Watson's \cite{W36} third order mock theta function $\nu(q)$.

\begin{Theorem}\label{tbsum}
If $|q|<1$, then
{\allowdisplaybreaks
\begin{align}\label{phibsum}
\phi(q)+\sum_{r=1}^{\infty}(-1;q^2)_rq^r&=\sum_{n=-\infty}^{\infty} \frac{q^{n^2}}{(-q^2;q^2)_n}=\frac{(-q,-q,q^2;q^2)_{\infty}}
{(q,-q^2;q^2)_{\infty}};\\
\nu(q)+\sum_{r=0}^{\infty}(-q;q^2)_rq^r&=\sum_{r=-\infty}^{\infty}\frac{q^{n^2+n}}{(-q;q^2)_{n+1}}
=2(-q^2,-q^2;q^2)_{\infty}(q^4;q^4)_{\infty};\numberthis \label{nubsum}\\
\psi(q)+\sum_{r=0}^{\infty}(q;q^2)_r(-1)^r&=\sum_{n=-\infty}^{\infty} \frac{q^{n^2}}{(q;q^2)_n}=\frac{(-q,-q,q^2;q^2)_{\infty}}
{2(q,-q^2;q^2)_{\infty}}. \numberthis \label{psibsum}
\end{align}}
\end{Theorem}

\begin{proof}
From \eqref{8baileyeqc12} (replace $q$ with $q^2$,
set $b=-z/t$, $a=-z$, and let $c\to \infty$),
\begin{equation*}
\sum_{r=-\infty}^{\infty}
\frac{(-z/t;q)_rt^r q^{r(r+1)/2}}{(tq;q)_r}
=\frac{(-t^2q^2/z,-zq,-q/z,q^2;q^2)_{\infty}}{(tq,-tq/z;q)_{\infty}},
\end{equation*}
and from \eqref{gzeq2} (with $s=-t$),
\begin{equation*}
\sum_{n=-\infty}^{\infty}\frac{z^n q^{n^2}}{(t^2q^2;q^2)_n}= \frac{(-tq/z;q)_{\infty}}{(-tq,-t^2q/z;q)_{\infty}}
\sum_{r=-\infty}^{\infty}
\frac{(-z/t;q)_r (t)^r q^{r(r+1)/2}}{(tq;q)_r}.
\end{equation*}
Together, these equations imply that
\begin{equation}
\sum_{n=-\infty}^{\infty}\frac{z^n q^{n^2}}{(t^2q^2;q^2)_n}=
\frac{(-zq,-q/z,q^2;q^2)_{\infty}}{(t^2q^2,-t^2q/z;q^2)_{\infty}}.
\end{equation}
The identity at \eqref{phibsum} is now immediate upon setting $z=1$ and $t^2=-1$, and that at
\eqref{psibsum} results similarly upon setting $z=1$ and $t^2=1/q$. The identity at \eqref{nubsum} follows upon setting $z=q$, $t^2=-q$, multiplying the resulting product by $1/(1+q)$, and finally performing some elementary $q$-product manipulations.
\qed
\end{proof}
Note that the convergence of the sum added to $\psi(q)$ on the left side of \eqref{psibsum} is in the Ces\`aro sense. Note also that comparison of the infinite products on the right sides of \eqref{phibsum} and \eqref{psibsum} yields the rather curious identity
\begin{equation}
\sum_{n=-\infty}^{\infty} \frac{q^{n^2}}{(-q^2;q^2)_n}=2\sum_{n=-\infty}^{\infty} \frac{q^{n^2}}{(q;q^2)_n},
\end{equation}
where, by the previous comment, convergence of the part of the bilateral series on the right consisting of terms of negative index is again in the Ces\`aro sense.

The summation formulae in the preceding theorem have some interesting implications. Firstly, they allow condition (0) above to be made explicit for some of the third order mock theta functions. We recall the recent result for $f(q)$  in \cite{FOR13}.

\begin{Theorem}\label{fort}(Folsom, Ono and Rhoades \cite{FOR13})
If $\zeta$ is a primitive even-order $2k$ root of unity, then, as $q$
approaches $\zeta$ radially within the unit disk, we have that
\begin{equation}\label{foreq}
\lim_{q\to \zeta}
(f(q)-(-1)^k b(q))= -4\sum_{n=0}^{k-1} (1+\zeta)^2(1+\zeta^2)^2\dots (1+\zeta^n)^2\zeta^{n+1}.
\end{equation}
\end{Theorem}
Here
\[
b(q)=\frac{(q;q)_{\infty}}{(-q;q)^2_{\infty}}.
\]
The infinite product representation of $b(q)$ was not stated in \cite{FOR13}, but was stated by Rhoades in  \cite{R13}. Note that Theorem \ref{fort} was also proved recently by Zudilin \cite{Z15}.

The following results are immediate upon rearranging the identities in Theorem \ref{tbsum}, and letting $q$ tend radially to the specified root of unity from within the unit circle, since the other series accompanying each of the mock theta functions in the bilateral sums terminates (the interchange of summation and limit in each of the corresponding series on the right  is justified by the absolute convergence of each of these series).
\begin{Corollary}\label{cmt3lim}
(i) If $\zeta$ is a primitive even-order $4k$ root of unity, then, as $q$
approaches $\zeta$ radially within the unit disk, we have that
\begin{equation}\label{philim}
\lim_{q\to \zeta}
\left(\phi(q)- \frac{(q^2,-q,-q;q^2)_{\infty}}{(-q^2,q;q^2)_{\infty}}\right)= -2\sum_{n=0}^{k-1} (1+\zeta^2)(1+\zeta^4)\dots (1+\zeta^{2n})\zeta^{n+1}.
\end{equation}
(ii) If $\zeta$ is a primitive even-order $4k+2$ root of unity, then, as $q$
approaches $\zeta$ radially within the unit disk, we have that
\begin{equation}\label{nulim}
\lim_{q\to \zeta}
\left(\nu(q)- 2(-q^2;q^2)^2_{\infty}(q^4;q^4)_{\infty}\right)= -\sum_{n=0}^{k} (1+\zeta)(1+\zeta^3)\dots (1+\zeta^{2n-1})\zeta^{n}.
\end{equation}
(iii) If $\zeta$ is a primitive odd-order $2k+1$ root of unity, then, as $q$
approaches $\zeta$ radially within the unit disk, we have that
\begin{equation}\label{psilim}
\lim_{q\to \zeta}
\left(\psi(q)- \frac{(q^2,-q,-q;q^2)_{\infty}}{2(-q^2,q;q^2)_{\infty}}\right)= -\sum_{n=0}^{k} (1-\zeta)(1-\zeta^3)\dots (1-\zeta^{2n-1})(-1)^{n}.
\end{equation}
\end{Corollary}
Remark: The results in Corollary \ref{cmt3lim} were also proved in \cite{BKLMR15}, using somewhat similar arguments, as were the results in Corollary \ref{cmt5lim}  below.

The second implication is that they imply some summation formulae for some of the bilateral series appearing in Theorem \ref{mt3bsum1}.

\begin{Corollary}\label{cmt3bsums}
If $|q|<1$, then
{\allowdisplaybreaks
\begin{align}\label{mt3bsum2eq1}
&\sum_{r=-\infty}^{\infty}
\frac{ q^{r^2+r}(2rq^{r}+1)}{(1+q^{r})^2}=
\frac{(q;q^2)^4_{\infty}(q;q)^4_{\infty}}{4}\\
&\sum_{r=-\infty}^{\infty}
\frac{ q^{2r^2+2r}(r+1)}{(1+q^{2r+1})^2}=
\frac{(-q^2;q^2)^2_{\infty}(q;q)^3_{\infty}(q^4;q^4)_{\infty}}{2(-q;q^2)_{\infty}}
+\frac{(-q^4,-q^{12},q^{16};q^{16})_{\infty}}{2}\numberthis \label{mt3bsum2eq2}
\end{align}}
\end{Corollary}
\begin{proof}
The first identity \eqref{mt3bsum2eq1} follows from combining the results at \eqref{phiq3id1} and \eqref{phibsum} and then replacing $q^2$ with $q$. The identity at \eqref{mt3bsum2eq2} follows directly from comparing the identities \eqref{nuq3id1} and \eqref{nubsum}.
\qed
\end{proof}

\section{  Mock theta functions of the fifth order}
\setcounter{equation}{0}

Ramanujan's fifth order mock theta functions are the following:
{\allowdisplaybreaks
\begin{align*}
f_0(q)&=\sum_{n=0}^{\infty} \frac{q^{n^2}}{(-q;q)_n},& f_1(q)&=\sum_{n=0}^{\infty} \frac{q^{n(n+1)}}{(-q;q)_n},\\
F_0(q)&=\sum_{n=0}^{\infty} \frac{q^{2n^2}}{(q;q^2)_n},& F_1(q)&=\sum_{n=0}^{\infty} \frac{q^{2n(n+1)}}{(q;q^2)_{n+1}},\\
\phi_0(q)&=\sum_{n=0}^{\infty} q^{n^2}(-q;q^2)_n, &\phi_1(q)&=\sum_{n=0}^{\infty} q^{(n+1)^2}(-q;q^2)_n,\\
\psi_0(q)&=\sum_{n=0}^{\infty} q^{(n+1)(n+2)/2}(-q;q)_n, &\psi_1(q)&=\sum_{n=0}^{\infty} q^{n(n+1)/2}(-q;q)_n,\\
\chi_0(q)&=\sum_{n=0}^{\infty} \frac{q^{n}(q;q)_n}{(q;q)_{2n}},& \chi_1(q)&=\sum_{n=0}^{\infty} \frac{q^{n}(q;q)_n}{(q;q)_{2n+1}}.
\end{align*}
}

Of interest here is the fact that that certain combinations of pairs of mock theta functions of order five may be expressed as single bilateral series, and hence in terms of theta products, as was described by Watson in section 7 of \cite{W37} (see also the forthcoming book \cite{McL16}, where  the proofs of these identities are possibly more transparent than those of Watson \cite{W37}). We state these identities directly in terms of $q$-products, rather than employing the Ramanujan functions $G(q)$ and $H(q)$, as Watson did.
\begin{Proposition}\label{28c1}
The following identities hold.
\begin{equation}\label{28c1eq1}
\sum_{n=-\infty}^{\infty}\frac{q^{r^2}}{(-q;q)_r}
=f_0(q)+2\psi_0(q)=4q\frac{(q^4,q^{16},q^{20};q^{20})_{\infty}}{(q^2;q^4)_{\infty}} +\frac{(q^2,q^{3},q^{5};q^{5})_{\infty}}{(-q;q)_{\infty}}.
\end{equation}
\begin{equation}\label{28c1eq2}
\sum_{n=-\infty}^{\infty}\frac{q^{r^2+r}}{(-q;q)_r}
=f_1(q)+2\psi_1(q)=4\frac{(q^8,q^{12},q^{20};q^{20})_{\infty}}{(q^2;q^4)_{\infty}} -\frac{(q,q^{4},q^{5};q^{5})_{\infty}}{(-q;q)_{\infty}}.
\end{equation}
\begin{multline}\label{28c1eq3}
\sum_{n=-\infty}^{\infty}\frac{q^{4r^2}}{(q^2;q^4)_r}\\
=F_0(q^2)+\phi_0(-q^2)-1=q\frac{(q^4,q^{16},q^{20};q^{20})_{\infty}}{(q^2;q^4)_{\infty}} +\frac{(q^2,q^{3},q^{5};q^{5})_{\infty}}{(-q;q)_{\infty}}.
\end{multline}
\begin{equation}\label{28c1eq4}
\sum_{n=-\infty}^{\infty}\frac{q^{4r^2+4r}}{(q^2;q^4)_{r+1}}
=F_1(q^2)-\frac{\phi_1(-q^2)}{q^2}=\frac{(q^8,q^{12},q^{20};q^{20})_{\infty}}{q(q^2;q^4)_{\infty}} -\frac{(q,q^{4},q^{5};q^{5})_{\infty}}{q(-q;q)_{\infty}}.
\end{equation}
\end{Proposition}

We note that these summation formulae may be rearranged and used to give explicit radial limits
for the difference of certain fifth order mock theta functions and certain corresponding theta functions,  as $q$ tends to certain roots of unity from within the unit circle,  in a manner similar to the result of Folsom, Ono and Rhoades \cite{FOR13} for the third order mock theta function $f(q)$ stated at \eqref{foreq} above, or to the results stated for the third order mock theta functions $\phi(q)$, $\nu(q)$ and $\psi(q)$ stated in Corollary \ref{cmt3lim}. For ease of notation, the statements for $F_0(q)$ and $F_1(q)$ are written in terms of $q^2$ instead of $q$.

{\allowdisplaybreaks
\begin{Corollary}\label{cmt5lim}
(i) If $\zeta$ is a primitive even-order $2k$ root of unity, then, as $q$
approaches $\zeta$ radially within the unit disk, we have that
\begin{multline}\label{f0lim}
\lim_{q\to \zeta}
\left(f_0(q)- \left[4q\frac{(q^4,q^{16},q^{20};q^{20})_{\infty}}{(q^2;q^4)_{\infty}} +\frac{(q^2,q^{3},q^{5};q^{5})_{\infty}}{(-q;q)_{\infty}}\right]\right)\\
= -2\sum_{n=0}^{k-1} (1+\zeta)(1+\zeta^2)\dots (1+\zeta^{n})\zeta^{(n+1)(n+2)/2}.
\end{multline}
(ii) If $\zeta$ is a primitive even-order $2k$ root of unity, then, as $q$
approaches $\zeta$ radially within the unit disk, we have that
\begin{multline}\label{f1lim}
\lim_{q\to \zeta}
\left(f_1(q)-  \left[4\frac{(q^8,q^{12},q^{20};q^{20})_{\infty}}{(q^2;q^4)_{\infty}} -\frac{(q,q^{4},q^{5};q^{5})_{\infty}}{(-q;q)_{\infty}}\right]\right)\\
= -2\sum_{n=0}^{k-1} (1+\zeta)(1+\zeta^2)\dots (1+\zeta^{n})\zeta^{n(n+1)/2}.
\end{multline}
(iii) If $\zeta$ is a primitive even-order $4k+2$ root of unity, then, as $q$
approaches $\zeta$ radially within the unit disk, we have that
\begin{multline}\label{F0lim}
\lim_{q\to \zeta}
\left(F_0(q^2)- \left[q\frac{(q^4,q^{16},q^{20};q^{20})_{\infty}}{(q^2;q^4)_{\infty}} +\frac{(q^2,q^{3},q^{5};q^{5})_{\infty}}{(-q;q)_{\infty}}\right]\right)\\
= -\sum_{n=1}^{k} (1-\zeta^2)(1-\zeta^6)\dots (1-\zeta^{4n-2})(-1)^n\zeta^{2n^2}.
\end{multline}
(iv) If $\zeta$ is a primitive even-order $4k+2$ root of unity, then, as $q$
approaches $\zeta$ radially within the unit disk, we have that
\begin{multline}\label{F1lim}
\lim_{q\to \zeta}
\left(F_1(q^2)- \left[\frac{(q^8,q^{12},q^{20};q^{20})_{\infty}}{q(q^2;q^4)_{\infty}} -\frac{(q,q^{4},q^{5};q^{5})_{\infty}}{q(-q;q)_{\infty}}\right]\right)\\
= -\sum_{n=0}^{k} (1-\zeta^2)(1-\zeta^6)\dots (1-\zeta^{4n-2})(-1)^n\zeta^{2n^2+4n}.
\end{multline}
\end{Corollary}
}

There are no known unilateral transformations for mock theta functions of the fifth order similar to those that exist for mock theta functions of the third order. However, there are bilateral transformations that may be applied to the bilateral series in Proposition \ref{28c1}.

Here we consider the series
\begin{equation}\label{g5def}
G_5^{*}(w,y,q)=\sum_{n=-\infty}^{\infty}\frac{w^nq^{n^2}}{(y;q)_n}.
\end{equation}
Identities for this function are not so plentiful as those for $G_3(s,t,q)$ and $G_3^{*}(s,t,q)$, but two such are given in the next theorem.
\begin{Theorem}\label{13tst}
Let $G_5^{*}(w,y,q)$ be as defined at \eqref{g5def}. Then
\begin{align}\label{g5eq1}
G_5^{*}(w,y,q)&=\frac{(y/w;q)_{\infty}}{(y;q)_{\infty}}
\sum_{r=-\infty}^{\infty}(wq/y;q)_r(-y)^rq^{r(r-1)/2},\\
&=\frac{(y/w;q)_{\infty}}{(wq,q/w;q)_{\infty}}
\sum_{r=-\infty}^{\infty}\frac{(1-w q^{2r})(wq/y;q)_r(-yw^2)^rq^{(5r^2-3r)/2}}
{(1-w)(y;q)_r}.\numberthis \label{g5eq2}
\end{align}
\end{Theorem}

\begin{proof}
In \eqref{q7.11eq1} (or  \eqref{q7.11eq2}), replace $z$ with $z/ac$, then let $a$, $c\to \infty$ and $d\to 0$. Then replace $z$ with $wq$ and $b$ with $y$, and \eqref{g5eq1} follows.

For \eqref{g5eq1}, let $d, e, f \to \infty$ in \eqref{8b2psi2}, and then set $a=w$ and $c= wq/y$.
\qed
\end{proof}
Remark: For $G_5^{*}(w,y,q)$ to represent a sum of fifth order mock theta function, it necessary to have $w=1$ or $w=q$, and in those cases \eqref{g5eq1} does not provide any non-trivial results (for $w=1$ the right side is just the series in reverse order).

The identity at \eqref{g5eq2} could be used to derive new expressions for the sums of fifth order mock theta functions found
in Corollary \ref{28c1}. However, we instead use it to derive four identities for bilateral series similar to those in Corollary \ref{cmt3bsums}.

\begin{Corollary}\label{28c2}
The following identities hold for $|q|<1$:
{\allowdisplaybreaks
\begin{align}\label{28c2eq1}
&\sum_{r=-\infty}^{\infty}(10r+1) q^{(5r^2+r)/2}
\\&
\phantom{asdasda}=\left(
\frac{4q(q^4,q^{16},q^{20};q^{20})_{\infty}}{(q^2;q^{4})_{\infty}}
+
\frac{(q^2,q^{3},q^{5};q^{5})_{\infty}}{(-q;q)_{\infty}}
\right)
\frac{(q;q)^2_{\infty}}{(-q;q)_{\infty}}
,\notag\\
&\sum_{r=-\infty}^{\infty}(10r+3) q^{(5r^2+3r)/2} \numberthis\label{28c2eq2}
\\&
\phantom{asdasda}=\left(
\frac{4(q^8,q^{12},q^{20};q^{20})_{\infty}}{(q^2;q^{4})_{\infty}}
-
\frac{(q,q^{4},q^{5};q^{5})_{\infty}}{(-q;q)_{\infty}}
\right)
\frac{(q;q)^2_{\infty}}{(-q;q)_{\infty}}
, \notag\\
&\sum_{r=-\infty}^{\infty}(5 r +1)(-1)^r q^{10r^2+4r} \numberthis \label{28c2eq3}
 \\
&\phantom{asdasda}=\bigg(
\frac{q(q^4,q^{16},q^{20};q^{20})_{\infty}}{(q^2;q^{4})_{\infty}}
+
\frac{(q^2,q^{3},q^{5};q^{5})_{\infty}}{(-q;q)_{\infty}}
\bigg)
\frac{(q^4;q^4)^2_{\infty}}{(q^2;q^4)_{\infty}}
, \notag\\
&\sum_{r=-\infty}^{\infty}(5r+2) (-1)^r q^{10r^2+8r}   \numberthis \label{28c2eq4}
\\&
\phantom{asdasda}=\left(
\frac{(q^8,q^{12},q^{20};q^{20})_{\infty}}{(q^2;q^{4})_{\infty}}
-
\frac{(q,q^{4},q^{5};q^{5})_{\infty}}{(-q;q)_{\infty}}
\right)
\frac{(q^4;q^4)^2_{\infty}}{(q^2;q^4)_{\infty}}
. \notag
\end{align}
}
\end{Corollary}
\begin{proof}
For \eqref{28c2eq1}, in \eqref{g5eq2} replace $w$ with $w^2$ , set $y=-wq$ and simplify the resulting right side to get
{\allowdisplaybreaks
\begin{align*}
&\sum_{n=-\infty}^{\infty}\frac{w^{2n}q^{n^2}}{(-wq;q)_n}
=\frac{(-q/w;q)_{\infty}}{(w^2q,q/w^2;q)_{\infty}}
\sum_{r=-\infty}^{\infty}\frac{(1-w q^{r})w^{5r}q^{(5r^2-r)/2}}
{(1-w)}\\
&=\frac{(-q/w;q)_{\infty}}{(w^2q,q^2/w;q)_{\infty}}\\
&\times
\left(
1+\sum_{r=1}^{\infty}\frac{(1-w q^{r})w^{5r}q^{(5r^2-r)/2}
+(1-w q^{-r})w^{-5r}q^{(5r^2+r)/2}}
{1-w}
\right)\\
&=\frac{(-q/w;q)_{\infty}}{(w^2q,q^2/w;q)_{\infty}}\\
&\times
\left(
1+\sum_{r=1}^{\infty}\frac{w^{-5r}q^{(5r^2-r)/2}}
{1-w}((1-w q^{r})w^{10r}+(1-w q^{-r})q^r)
\right)\\
&=\frac{(-q/w;q)_{\infty}}{(w^2q,q^2/w;q)_{\infty}}\\
&\times
\left(
1+\sum_{r=1}^{\infty}w^{-5r}q^{(5r^2-r)/2}
\left(-w\frac{1-w^{10r-1}}{1-w}+q^r \frac{1-w^{10r+1}}{1-w}\right)
\right).
\end{align*}
}
Now let $w \to 1$, noting that the left side above tends to the left side of \eqref{28c1eq1}, and hence to the right side of \eqref{28c1eq1}. After using L'Hospital's rule on the terms in the last series on the right side, this series becomes
{\allowdisplaybreaks
\begin{align*}
1&+\sum_{r=1}^{\infty}q^{(5r^2-r)/2}
\left(-(10r-1)+q^r(10r+1)\right)\\
&=\sum_{r=-\infty}^{\infty} 10r q^{(5r^2+r)/2}+\sum_{r=-\infty}^{\infty}  q^{(5r^2-r)/2}.
\end{align*}
}
The result now follows.

To obtain \eqref{28c2eq2}, in \eqref{g5eq2} replace $w$ with $w^2q$ , set $y=-wq$ and simplify the resulting right side to get
{\allowdisplaybreaks
\begin{align*}
&\sum_{n=-\infty}^{\infty}\frac{w^{2n}q^{n^2+n}}{(-wq;q)_n}
=\frac{(-1/w;q)_{\infty}}{(w^2q^2,1/w^2;q)_{\infty}}
\sum_{r=-\infty}^{\infty}\frac{(1-w^2 q^{2r+1})w^{5r}q^{(5r^2+3r)/2}}
{(1-w^2q)}\\
&=\frac{(-1/w;q)_{\infty}(-w^2)}{(w^2q,q/w^2;q)_{\infty}}\times\\
&
\sum_{r=0}^{\infty}\frac{(1-w^2 q^{2r+1})w^{5r}q^{(5r^2+3r)/2}}
{(1-w^2)}
+
\sum_{r=0}^{\infty}\frac{(1-w^2 q^{-2r-1})w^{-5r-5}q^{(5r^2+7r+2)/2}}
{(1-w^2)}\\
&=\frac{(-1/w;q)_{\infty}(-w^2)}{(w^2q,q/w^2;q)_{\infty}}\times\\
&
\sum_{r=0}^{\infty}w^{5r}q^{(5r^2+3r)/2}
\left(\frac{(1-w^{-10r-3})}
{(1-w^2)}
-w^2q^{2r+1}
\frac{(1-w^{-10r-7})}
{(1-w^2)}\right),
\end{align*}
}
where the second series in the second right side came from taking the terms of negative index in the series on the first right side, and replacing $r$ with $-r-1$.
The identity at \eqref{28c2eq2} now follows as previously upon letting $w \to 1$, this time noting that the left side tends to \eqref{28c1eq2}.

For \eqref{28c2eq3} and \eqref{28c2eq4}, in \eqref{g5eq2} replace $(w,y,q)$ with $(w^2, wq^2, q^4)$ and $(w^2q^4, wq^6,$ $ q^4)$, respectively (in the case of \eqref{28c2eq4}, after making the replacements in \eqref{g5eq2}, multiply both sides by $1/(1-wq^2)$). The details are similar to those in the proofs of \eqref{28c2eq1} and \eqref{28c2eq2}, and are omitted.
\qed
\end{proof}

\section{  Mock theta functions of the sixth order}
\setcounter{equation}{0}

The sixth order mock theta functions which concern us here are
{\allowdisplaybreaks
\begin{align}
&\phi(q)=\sum_{n=0}^{\infty}\frac{ (-1)^n q^{n^2}(q;q^2)_n}{(-q;q)_{2n}},
&\psi(q)=\sum_{n=0}^{\infty} \frac{(-1)^n q^{(n+1)^2}(q;q^2)_n}{(-q;q)_{2n+1}},\notag\\
&\rho(q)=\sum_{n=0}^{\infty} \frac{q^{n(n+1)/2}(-q;q)_n}{(q;q^2)_{n+1}},
&\sigma(q)=\sum_{n=0}^{\infty} \frac{q^{(n+1)(n+2)/2}(-q;q)_n}{(q;q^2)_{n+1}},\notag\\
&\lambda(q)=\sum_{n=0}^{\infty} \frac{(-1)^n q^{n}(q;q^2)_{n}}{(-q;q)_{n}},
&\mu(q)=\sum_{n=0}^{\infty} \frac{(-1)^n(q;q^2)_n}{(-q;q)_n},\notag\\
&\phi_{-}(q)=\sum_{n=1}^{\infty} \frac{q^{n}(-q;q)_{2n-1}}{(q;q^2)_n},
&\psi_{-}(q)=\sum_{n=1}^{\infty} \frac{q^{n}(-q;q)_{2n-2}}{(q;q^2)_n}.\notag
\end{align}
}

The series stated for $\mu(q)$ \label{mupage} does not converge, but  the sequence of even-indexed partial sums and the sequence of odd-indexed partial sums do converge, and $\mu(q)$ is defined to be the average of these two values.

Andrews and Hickerson \cite{AH91} proved a number of identities for the sixth order mock theta functions  stated by Ramanujan in the Lost Notebook \cite{R88}.
Berndt and Chan \cite{BC07} proved a number of similar identities. The proofs in both of these papers were quite involved, employing both Bailey pairs and the constant term method,  and simpler proofs were later given by Lovejoy \cite{L10}, for four of the identities proved by Andrews and Hickerson \cite{AH91}.
These four identities are listed in the following theorem.
\begin{Theorem}\label{mt6t1}
The following identities hold for $|q|<1$.
{\allowdisplaybreaks
\begin{align}\label{mt6t1eq1}
q^{-1}\psi(q^2)+\rho(q)&=(-q;q^2)^2_{\infty}(-q,-q^5,q^6;q^6)_{\infty},\\
\phi(q^2)+2\sigma(q)&=(-q;q^2)^2_{\infty}(-q^3,-q^3,q^6;q^6)_{\infty},\numberthis \label{mt6t1eq2}\\
2\phi(q^2)-2\mu(-q)&=(-q;q^2)^2_{\infty}(-q^3,-q^3,q^6;q^6)_{\infty},\numberthis \label{mt6t1eq3}\\
2q^{-1}\psi(q^2)+\lambda(-q)&=(-q;q^2)^2_{\infty}(-q,-q^5,q^6;q^6)_{\infty}, \numberthis \label{mt6t1eq4}
\end{align}
}
\end{Theorem}

To maintain uniformity, we show that, as with the third order- and fifth order mock theta functions, the bilateral transformations of Bailey at \eqref{q7.11eq1}, \eqref{q7.11eq2} and \eqref{8b2psi2} may be used to express sums of sixth order mock theta functions as theta functions, and that these identities in turn may likewise be used to examine the limiting behavior of some of these sixth order mock theta functions as $q$ tends to certain classes of roots of unity from within the unit circle. As above, we begin by stating a number of general bilateral transformations.

\begin{Theorem}\label{mt6t2}
(i) If $|q|, |bd/azq|<1$, then
{\allowdisplaybreaks
\begin{multline}\label{28m6eq1}
G_6(a,b,d,z,q):=\sum_{n=-\infty}^{\infty}\frac{(a;q^2)_r z^r q^{r^2}}{(b,d;q^2)_r}
=\frac{(-zq,-qb/az,-qd/az;q^2)_{\infty}}{(b,d,q^2/a;q^2)_{\infty}}\\
\times\sum_{n=-\infty}^{\infty}\frac{(-azq/b,-azq/d;q^2)_r}{(-zq;q^2)_r}
\left(
\frac{-bd}{azq}
\right)^r.
\end{multline}
(ii) If $|q|, |bd/azq|, |b/a|<1$, then
\begin{multline}\label{28m6eq2}
\sum_{n=-\infty}^{\infty}\frac{(a;q^2)_r z^r q^{r^2}}{(b,d;q^2)_r}
=\frac{(b/a,-qb/az;q^2)_{\infty}}{(b,-bd/azq;q^2)_{\infty}}\\
\times\sum_{n=-\infty}^{\infty}\frac{(-azq/b,a;q^2)_r}{(d;q^2)_r}
\left(
\frac{b}{a}
\right)^r.
\end{multline}
(iii) If $|q|, |bd/azq|<1$, then
\begin{multline}\label{28m6eq3}
\sum_{n=-\infty}^{\infty}\frac{(a;q^2)_r z^r q^{r^2}}{(b,d;q^2)_r}
=\frac{(-bq/az,-dq/az,-qz;q^2)_{\infty}}{(-bd/aqz,-q^3/az,-aqz;q^2)_{\infty}}\\
\times\sum_{n=-\infty}^{\infty}\frac{(1+azq^{4r-1})(a,-azq/b,-azq/d;q^2)_r\left(
bdz
\right)^rq^{3r^2-4r}}
{(1+az/q)(b,d,-zq;q^2)_r}.
\end{multline}
}
\end{Theorem}
\begin{proof}
For \eqref{28m6eq1} and \eqref{28m6eq2}, replace $q$ with $q^2$, $z$ with $-zq/c$ and then let $c\to \infty$ in \eqref{q7.11eq1} and \eqref{q7.11eq2}, respectively.

The identity \eqref{28m6eq3} is a consequence of replacing $q$ with $q^2$ in \eqref{8b2psi2}, and then replacing $c$ with $aq^2/b$ and $d$ with $aq^2/d$, letting $f \to \infty$, and replacing, in turn, $a$ with $-ze/q$ and finally $e$ with $a$.
\qed
\end{proof}

The sums of various pairs of sixth order mock theta functions may be expressed in terms of $G_6(a,b,d,z,q)$, and the above theorem may be used to derive  some alternative expressions for these sums.

\begin{Corollary}\label{28cgzpm}
The following identities hold for $|q|<1$.
{\allowdisplaybreaks
\begin{align}\label{28cgzpmeq0}
&4\sigma(q)+2\mu(q)=\frac{(q;q^2)_{\infty}}{(q;q)^2_{\infty}}
\sum_{r=-\infty}^{\infty}(-1)^r(6r+1)q^{r(3r+1)/2},\\
&\phi(q)+2\phi_{-}(q)=\frac{(-q;q)_{\infty}}{(q^2;q^4)_{\infty}}\numberthis \label{28cgzpmeq1}\\
&\phantom{}\times \big[
2(-q^2;q^4)^2_{\infty}(-q^6,-q^6,q^{12};q^{12})_{\infty}
-(q^2;q^4)^2_{\infty}(q^6,q^6,q^{12};q^{12})_{\infty}
\big],\notag\\
&\psi(q)+2\psi_{-}(q)=\frac{3q(-q;q)_{\infty}(q^6,q^6,q^6;q^6)_{\infty}}
{(q^2;q^2)^2_{\infty}}\numberthis \label{28cgzpmeq2}\\
&2\rho(q)+\lambda(q)=\frac{3(q;q^2)_{\infty}(q^3,q^3,q^3;q^3)_{\infty}}
{(q;q)^2_{\infty}}\numberthis \label{28cgzpmeq3}
\end{align}
}
\end{Corollary}

\begin{proof}
In \eqref{28m6eq3} replace $z$ with $-zq^3$ and set $a=zq^2$, $b=zq^3$ and $d=-zq^3$ to get
{\allowdisplaybreaks
\begin{align*}
&\sum_{r=-\infty}^{\infty}
\frac{(zq^2;q^2)_r (-z)^r q^{r^2+3r}}
{(zq^3,-zq^3;q^2)_r}=
\frac{(1/qz,-1/qz,zq^4;q^2)_{\infty}}
{(-1,1/q^2z^2,q^6z^2;q^2)_{\infty}}\\
&\times
\sum_{r=-\infty}^{\infty}
\frac{q^{3 r^2+5 r} z^{3 r} \left(1-z^2
   q^{4 r+4}\right) \left(q^2
   z;q^2\right){}_r}{\left(1-q^4 z^2\right)
   \left(q^4 z;q^2\right){}_r},\\
&=
\frac{(q^2/z^2;q^4)_{\infty}(zq^4;q^2)_{\infty}z}
{2(1-1/z)(-q^2,q^2/z^2,q^6z^2;q^2)_{\infty}}
\sum_{r=-\infty}^{\infty}
\frac{ q^{r (3 r+5)} z^{3 r} \left(1+z
   q^{2 r+2}\right)}{(1+q^2 z)(1+z)}
\end{align*}
}
Now multiply both sides by $q^2/(1-z^2q^2)$ and let $z\to -1$, noting that the left side tends to $\sigma(q^2)+\mu(q^2)/2$ using the definitions above and \eqref{1def4}. On the right side replace $r$ with $r-1$ and rewrite the resulting series  as
\begin{multline*}
\sum_{r=-\infty}^{\infty}
\frac{(1+z q^{2r})z^{3r}q^{3r^2-r}}
{1+z}
=\sum_{r=-\infty}^{\infty}
\frac{z^{3r}q^{3r^2-r}+z^{3r+1}q^{3r^2+r}}
{1+z}\\
=\sum_{r=-\infty}^{\infty}
\frac{z^{3r}q^{3r^2-r}+z^{-3r+1}q^{3r^2-r}}
{1+z}
=\sum_{r=-\infty}^{\infty}z^{3r}q^{3r^2-r}
\frac{1+z^{-6r+1}}
{1+z},
\end{multline*}
where the second equality follows from reversing the order of summation for the second terms in the sum. Now let $z\to -1$ to arrive at
\[
\sigma(q^2)+\frac{\mu(q^2)}{2}=\frac{(q^2;q^4)_{\infty}}{4(q^2;q^2)^2_{\infty}}
\sum_{r=-\infty}^{\infty}(6r+1)(-1)^r q^{3r^2+r}.
\]
The identity at \eqref{28cgzpmeq0} now follows upon multiplying this last identity by 4 and replacing $q$ with $q^{1/2}$.

Note that the expression for $4\sigma(q)+2\mu(q)$ deriving from \eqref{mt6t1eq2} and \eqref{28cgzpmeq0}  together with \eqref{mt6t1eq2} imply that
{\allowdisplaybreaks
\begin{multline}\label{6r+1sum}
\sum_{r=-\infty}^{\infty}(6r+1)(-1)^r q^{(3r^2+r)/2}
=\frac{(q;q)^2_{\infty}}{(q;q^2)_{\infty}}\\
\times \left [
2(-q;q^2)^2_{\infty}(-q^3,-q^3,q^6;q^6)_{\infty}
-(q;q^2)^2_{\infty}(q^3,q^3,q^6;q^6)_{\infty}
\right].
\end{multline}
}

Remark: It may be of interest to compare the identity above with that of Fine \cite[p.83]{F88}:
\begin{equation}\label{fineid}
\sum_{r=-\infty}^{\infty} (6r+1)q^{r(3r+1)/2}=(q;q)_{\infty}^3(q;q^2)_{\infty}^2.
\end{equation}

For \eqref{28cgzpmeq1}, set $a=-zq$, $b=zq$ and $d=zq^2$ in \eqref{28m6eq3} to get, after simplifying the right side
\[
\sum_{r=-\infty}^{\infty}
\frac{(-zq;q^2)_r z^r q^{r^2}}
{(zq,zq^2;q^2)_r}=
\frac{(q/z,q^2/z,-zq;q^2)_{\infty}}
{(q,q^2/z^2,q^2z^2;q^2)_{\infty}}
\sum_{r=-\infty}^{\infty}
\frac{(1+z q^{2r})z^{3r}q^{3r^2-r}}
{1+z},
\]
Let $z\to -1$ on the left side to get, once again using the definitions above and \eqref{1def4},
\begin{align*}
\phi(q)+2\phi_{-}(q)=\frac{(-q;q)_{\infty}}{(q^2;q^2)^2_{\infty}}
\sum_{r=-\infty}^{\infty}(6r+1)(-1)^r q^{3r^2+r}.
\end{align*}
An application of \eqref{6r+1sum}, with $q$ replaced with $q^2$, gives the result.

Similarly, for \eqref{28cgzpmeq2}, replace $z$ with $zq^2$, $a$ with $-zq$, $b$ with $zq^2$ and $d$ with $zq^3$ in \eqref{28m6eq3} to get, after once again simplifying the right side, that
\begin{align*}
\sum_{r=-\infty}^{\infty}
\frac{(-zq;q^2)_r z^r q^{r^2+2r}}
{(zq^2,zq^3;q^2)_r}&=
\frac{(1/z,q/z,-zq^3;q^2)_{\infty}}
{(q,1/z^2,q^4z^2;q^2)_{\infty}}
\sum_{r=-\infty}^{\infty}
z^{3r}q^{3r^2+3r} \\
&=
\frac{
(1/z,q/z,-zq^3;q^2)_{\infty}(-q^6z^3,-1/z^3,q^6;q^6)_{\infty }
}
{
(q,1/z^2,q^4z^2;q^2)_{\infty}
},
\end{align*}
where the Jacobi triple product identity \eqref{6t1eq} has been used at the last step. The result now follows after multiplying both sides by $q/(1+q)$ and letting $z\to -1$ as before.

The details of the proof of \eqref{28cgzpmeq3} are omitted. Briefly, replace $z$ with $zq$, $a$ with $-zq^2$, $b$ with $zq^3$ and $d$ with $-zq^3$ in \eqref{28m6eq3}, simplify and sum the right side using the Jacobi triple product identity \eqref{6t1eq}, let $z\to 1$, multiply both sides by $2/(1-q^2)$, and finally replace $q$ with $q^{1/2}$.
\qed
\end{proof}

Remark: Choi \cite[p. 370]{Ch11} also gave expressions for each of the sums of sixth order mock theta functions in Corollary \ref{28cgzpm}, but with different combinations of theta functions on the right sides.
Yet another version of \eqref{28cgzpmeq1} was stated by Ramanujan \cite[p. 6 and p. 16]{R88} (see also \cite[p. 1740]{ChK12}). Different proofs of \eqref{28cgzpmeq2} and \eqref{28cgzpmeq3} were given by Choi and Kim \cite[Theorem 1.4, p. 1742]{ChK12}. The identities in Corollary \ref{28cgzpm} also follow from expressions  for the sixth order mock theta functions in terms of the function $m(x,q,z)$ (see \eqref{mxqzeq}) proved by Hickerson and  Mortenson in \cite{HM14},
and known results about $m(x,q,z)$.

We note that the identities in  Corollary \ref{28cgzpm} may be used to describe the asymptotic behavior of each of the two sixth order mock theta functions on the left side of each identity, at particular classes of roots of unity. For example, \eqref{28cgzpmeq0} may be used in conjunction with \eqref{6r+1sum}, to make condition (0) in the interpretation by Andrews and Hickerson  of a mock theta function explicit for both $\sigma(q)$ at primitive roots of unity of odd order, and for $\mu(q)$ at primitive roots of unity of even order.
We state the result for just one of each pair of mock theta functions, and leave the result for the other mock theta function of each pair to the reader.

{\allowdisplaybreaks
\begin{Corollary}\label{cmt6lim}
(i) If $\zeta$ is a primitive odd-order $2k+1$ root of unity, then, as $q$
approaches $\zeta$ radially within the unit disk, we have that
\begin{multline}\label{sig6lim}
\lim_{q\to \zeta}
\left(\!\sigma(q)- \frac{1}{4}\left[2(-q;q^2)^2_{\infty}(-q^3,-q^{3},q^{6};q^{6})_{\infty} -(q;q^2)^2_{\infty}(q^3,q^{3},q^{6};q^{6})_{\infty} \right]\!\right)\\
= -\frac{1}{2}\sum_{n=0}^{k} \frac{(1-\zeta)(1-\zeta^3)\dots (1-\zeta^{2n-1})}
{(1+\zeta)(1+\zeta^2)\dots (1+\zeta^{n})}
(-1)^{n}.
\end{multline}
(ii) If $\zeta$ is a primitive even-order $2k$ root of unity, then, as $q$
approaches $\zeta$ radially within the unit disk, we have that
\begin{multline}\label{phi6lim}
\lim_{q\to \zeta}
\bigg(\phi(q)-  \frac{(-q;q)_{\infty}}{(q^2;q^4)_{\infty}}\\
\times \bigg[
2(-q^2;q^4)^2_{\infty}(-q^6,-q^6,q^{12};q^{12})_{\infty}
-(q^2;q^4)^2_{\infty}(q^6,q^6,q^{12};q^{12})_{\infty}
\bigg]\bigg)\\
= -2\sum_{n=1}^{k} \frac{(1+\zeta)(1+\zeta^2)\dots (1+\zeta^{2n-1})}
{(1-\zeta)(1-\zeta^3)\dots (1-\zeta^{2n-1})}\zeta^{n}.
\end{multline}
(iii) If $\zeta$ is a primitive even-order $2k$ root of unity, then, as $q$
approaches $\zeta$ radially within the unit disk, we have that
\begin{multline}\label{psi6lim}
\lim_{q\to \zeta}
\left(\psi(q)- \frac{3q(-q;q)_{\infty}(q^6,q^6,q^6;q^6)_{\infty}}
{(q^2;q^2)^2_{\infty}}\right)\\
= -2\sum_{n=1}^{k} \frac{(1+\zeta)(1+\zeta^2)\dots (1+\zeta^{2n-2})}
{(1-\zeta)(1-\zeta^3)\dots (1-\zeta^{2n-1})}\zeta^{n}.
\end{multline}
(iv) If $\zeta$ is a primitive odd-order $2k+1$ root of unity, then, as $q$
approaches $\zeta$ radially within the unit disk, we have that
\begin{multline}\label{rho6lim}
\lim_{q\to \zeta}
\left(\rho(q)- \frac{3(q;q^2)_{\infty}(q^3,q^3,q^3;q^3)_{\infty}}
{2(q;q)^2_{\infty}}\right)\\
= -\frac{1}{2}\sum_{n=0}^{k} \frac{(1-\zeta)(1-\zeta^3)\dots (1-\zeta^{2n-1})}
{(1+\zeta)(1+\zeta^2)\dots (1+\zeta^{n})}(-\zeta)^{n}.
\end{multline}
\end{Corollary}
}

Before considering eighth order mock theta functions, we compare the results in the present paper with those implied by an identity of Mortenson  (\cite[Eq. (6.10)]{M15}):
\begin{multline}\label{morteq}
\sum_{n=0}^{\infty}\frac{q^{n(n+1)/2}(-q;q)_n}{(x;q)_{n+1}(q/x)_{n+1}}
+\sum_{n=0}^{\infty}\frac{1}{2}\frac{q^n(q/x;q)_n(x;q)_n}{(-q;q)_n}=-\frac{j(x;q)}{2J_2}g_3(-x;q)\\
+\frac{J_2^3}{J_{1,2}j(x^2;q^2)}+\frac{1}{2x}\frac{J_{2}^{10}j(-x^2;q^2)}{J_{1}^4J_{4}^4j(x^2;q^2)j(-qx^2;q^2)}
-\frac{1}{2x}\frac{J_{2,4}^2j(x;q)}{j(-x;q)j(-qx^2;q^2)},
\end{multline}
where  the first series is $g_2(x,q)$, the universal mock theta function of Gordon and McIntosh \cite[Eq. (4.11)]{GM12}, $g_3(x;q)$ is as defined at \eqref{g3xq}, and
\begin{align*}
&j(x;q):=(x,q/x,q;q)_{\infty},  &&J_{a,m}:=j(q^a;q^m),&\\
& \bar{J}_{a,m}:=j(-q^a;q^m),&& J_m:=J_{m,3m}=(q^m;q^m)_{\infty}.&
\end{align*}
As Mortenson indicated in \cite{M15}, if a mock theta function is expressible in terms of $g_2(x,q)$ and combinations of infinite products, then it may be possible to derive a radial limits result for certain classes of roots of unity, and indeed Mortenson derives such results for a second order mock theta function and one of tenth order, and states that there are many other cases where \eqref{morteq} may be applied.

As one way of deriving explicit radial limits, one might  hope, after substituting for $g_2(x,q)$ in \eqref{morteq}, so that this expression now contains a mock theta function, that there is then a class of roots of unity such that as $q$ approaches  one of these roots of unity, say $\zeta$, from within the unit circle, the mock theta function becomes unbounded, the term involving $g_3(-x;q)$ vanishes, and the second series on the left terminates. In this case \eqref{morteq} may then be rearranged to give an identity of the form
\[
\lim_{q\to \zeta}\,\,(\text{mock theta function} - \text{theta function})\,= \,\text{finite $q$-series in $\zeta$,}
\]
which is the typical form of a radial limits result.
For example, if $q$ is replaced with $q^6$ and $x$ with $q^3$ in \eqref{morteq} and the second identity at \cite[Eq. (5.10)]{GM12}, namely
\begin{equation}\label{psiq4eq}
\psi(q^4)=\frac{q^3J_2^2J_4J_24^2}{J_1J_3J_8^2}-q^3g_2(q^3,q^6)
\end{equation}
is used to substitute for $g_2(q^3,q^6)$, then after some $q$-product manipulation, we get
\begin{multline}\label{psiq4eq2}
\psi(q^4)+\frac{q^3J_{12}^5}{J_6^4}
+
\frac{J_{12}^{17}}{4J_6^8J_{24}^8}
-
\frac{J_3^4J_{12}^7}{4J_6^6J_{24}^4}
-
\frac{q^3J_2^2J_4J_{24}^2}{J_1J_3J_8^2}
-
\frac{q^3J_3^2}{2J_6J_{12}}
\sum_{n=0}^{\infty}\frac{q^{6n(n+1)}}{(-q^3;q^6)^2_{n+1}}\\
=\sum_{n=0}^{\infty}\frac{q^{6n+3}(q^3,q^3;q^6)_{n}}{(-q^6;q^6)_{n}}.
\end{multline}
If $q\to \zeta$, where $\zeta$ is a primitive even-order root of unity, then both the series on the right of \eqref{psiq4eq2} and the last term on the left become unbounded, and there is no radial limit. Unfortunately for producing explicit radial limits, when $q\to \zeta$, where $\zeta$ is a primitive \emph{odd}-order root of unity, while the series on the right of \eqref{psiq4eq2} terminates, and the last term on the left vanishes, the series for $\psi(q^4)$ also terminates. After eliminating terms that vanish when $q\to \zeta$, where $\zeta$ is a primitive $2k+1$-th root of unity, one gets that
\begin{multline}\label{psiq4eq3}
\lim_{q\to \zeta}\left(\frac{q^3J_{12}^5}{J_6^4}
+
\frac{J_{12}^{17}}{4J_6^8J_{24}^8}
-
\frac{q^3J_2^2J_4J_{24}^2}{J_1J_3J_8^2}\right)\\
=
\sum_{r=0}^{k}\frac{\zeta^{6r+3}(\zeta^3,\zeta^3;\zeta^6)_{r}}
{2(-\zeta^6;\zeta^6)_{r}}
-
\sum _{r=0}^{k} \frac{(-1)^r \zeta ^{4 (r+1)^2}
   \left(\zeta ^4;\zeta ^8\right){}_r}{\left(-\zeta ^4;\zeta
   ^4\right){}_{2 r+1}}.
\end{multline}
Curiously, what experiment suggests is that each side is identically zero, and in particular, that if $\zeta$ is a primitive $2k+1$-th root of unity, then
\begin{equation}\label{psiq4eq4}
\sum_{r=0}^{k}\frac{\zeta^{6r+3}(\zeta^3,\zeta^3;\zeta^6)_{r}}
{2(-\zeta^6;\zeta^6)_{r}}
=
\sum _{r=0}^{k} \frac{(-1)^r \zeta ^{4 (r+1)^2}
   \left(\zeta ^4;\zeta ^8\right){}_r}{\left(-\zeta ^4;\zeta
   ^4\right){}_{2 r+1}} \,\,\,\,(=\psi(\zeta^4))
\end{equation}
holds for all $k\geq 0$. Note that equality in \eqref{psiq4eq4} does not hold if $\zeta$ is replaced with an arbitrary value of $q$  inside the unit circle, since the difference is a non-zero function of $q$:
\begin{multline*}
\sum_{r=0}^{k}\frac{q^{6r+3}(q^3,q^3;q^6)_{r}}
{2(-q^6;q^6)_{r}}
-
\sum _{r=0}^{k} \frac{(-1)^r q ^{4 (r+1)^2}
   \left(q ^4;q ^8\right){}_r}{\left(-q ^4;q
   ^4\right){}_{2 r+1}}\\
   = \frac{q^3}{2}-q^4+q^8+\frac{q^9}{2}-2
   q^{12}+\frac{q^{15}}{2}+2 q^{16}-3
   q^{20}\dots .
\end{multline*}
Note also that each of the three individual terms on the left side of \eqref{psiq4eq3} diverges to $\infty$ as $q\to\zeta$, even though the combination converges to zero.
 Of course \eqref{psiq4eq4} will hold if
\[
\frac{q^3J_{12}^5}{J_6^4}
+
\frac{J_{12}^{17}}{4J_6^8J_{24}^8}
-
\frac{q^3J_2^2J_4J_{24}^2}{J_1J_3J_8^2}
= (q;q^2)_{\infty} \theta (q),
\]
where $\theta (q)$ is a function of $q$ that remains bounded as $q$ approaches any primitive odd-order root of unity from within the unit circle. We have not attempted to prove this, nor \eqref{psiq4eq4}.

If Ramanujan's identity (see \cite[Eq. (5.8)]{GM12})
\[
2q^{-1}\psi(q^2)+\lambda(-q)=(-q;q^2)_{\infty}^2j(-q,q^6)
\]
with $q$ replaced with $-q^2$ is used to replace $\psi(q^4)$ in \eqref{psiq4eq2}, and a similar analysis of radial limits is attempted, what experiment also appears to indicate is that
\begin{equation}\label{psiq4eq5}
\sum_{r=0}^{k}\frac{\zeta^{6r+1}(\zeta^3,\zeta^3;\zeta^6)_{r}}
{(-\zeta^6;\zeta^6)_{r}}
=
\sum _{r=0}^{k} \frac{(-1)^r \zeta ^{2r}
   \left(\zeta ^2;\zeta ^4\right){}_r}{\left(-\zeta ^2;\zeta
   ^2\right){}_{r}} \,\,\,\,(=\lambda(\zeta^2))
\end{equation}
holds, where $\zeta$ is a $2k+1$-th primitive root of unity (and again \eqref{psiq4eq5} does not hold if $\zeta$ is replaced with an arbitrary $q$ inside the unit circle). Similar results may be obtained from \eqref{morteq} for other sixth order mock theta functions.

\section{  Mock theta functions of the eighth order}
\setcounter{equation}{0}

We consider four of the eight mock theta functions of order eight   introduced by Gordon and McIntosh \cite{GM00}:
{\allowdisplaybreaks
\begin{align}\label{mt8defs1}
&S_0(q)=\sum_{n=0}^{\infty} \frac{q^{n^2}(-q;q^2)_n}{(-q^2;q^2)_{n}},
&S_1(q)=\sum_{n=0}^{\infty} \frac{q^{n(n+2)}(-q;q^2)_n}{(-q^2;q^2)_{n}},\\
&T_0(q)=\sum_{n=0}^{\infty}\frac{q^{(n+1)(n+2)}(-q^2;q^2)_n}{(-q;q^2)_{n+1}},
&T_1(q)=\sum_{n=0}^{\infty} \frac{q^{n(n+1)}(-q^2;q^2)_n}{(-q;q^2)_{n+1}}.
\numberthis \label{mt8defs2}
\end{align}
}

As with mock theta functions of other orders, certain sums of eighth order mock theta functions may be written as single bilateral series, and it is a straightforward consequence of the definitions and \eqref{1def4} that
\begin{align}\label{stbilatsumeq1}
S_0(q)+2T_0(q)&=\sum_{r=-\infty}^{\infty}
\frac{(-q;q^2)_rq^{r^2}}{(-q^2;q^2)_r},\\
S_1(q)+2T_1(q)&=\sum_{r=-\infty}^{\infty}
\frac{(-q;q^2)_rq^{r^2+2r}}{(-q^2;q^2)_r}.\numberthis \label{stbilatsumeq2}
\end{align}
In fact, following the method of the authors in \cite{GM00}, it will be shown that each of these sums has an expression in terms of infinite products. We include the proof here since in \cite{GM00} the authors  omitted the final step of explicitly stating the form of the infinite products (although these expressions were stated by them in \cite[Eq. (5.12)]{GM12}, and these expressions with further details of the proof were given by them in \cite[Section 4]{GM03}). As with the identities in Corollary \ref{28cgzpm}, the identities in Theorem \ref{st01t1} may also be shown to follow from identities proved by Hickerson and  Mortenson in \cite{HM14}.
\begin{Theorem}\label{st01t1}
If $|q|<1$, then
\begin{align}\label{stbilatsumprodeq1}
S_0(q^2)+2T_0(q^2)&=
\frac{(q^2;q^2)_{\infty}\left[(q;q^2)^3_{\infty}+  (-q;q^2)^3_{\infty}\right]}{2(-q^2;q^2)_{\infty}}
,\\
S_1(q^2)+2T_1(q^2)&=\frac{(q^2;q^2)_{\infty}\left[(-q;q^2)^3_{\infty}-  (q;q^2)^3_{\infty}\right]}{2q(-q^2;q^2)_{\infty}}
.\numberthis \label{stbilatsumprodeq2}
\end{align}
\end{Theorem}
\begin{proof}
Let $R_0(q)$ and $R_1(q)$ denote the series on the right side of \eqref{stbilatsumeq1} and \eqref{stbilatsumeq2}, respectively.
Next, in \eqref{8baileyeq1} replace $q$ with $q^2$, set $a=q$, $b=iq$ and $c=-iq$, and then let $d,e \to \infty$. This leads to
\begin{multline*}
\frac{R_0(q^2)-qR_1(q^2)}{1-q}=
\sum_{r=-\infty}^{\infty} \frac{(1-q^{4r+1})(-q^2;q^4)_rq^{2r^2}}{(1-q)(-q^4;q^4)_r}\\
=
\frac{(q^3,q,q^2,q;q^2)_{\infty}}{(iq^2,-iq^2,iq,-iq;q^2)_{\infty}}
=
\frac{(q^2;q^2)_{\infty}(q;q^2)^3_{\infty}}{(1-q)(-q^2;q^2)_{\infty}}.
\end{multline*}
Multiply through by $1-q$ and then replace $q$ with $-q$ to get an expression for $R_0(q^2)+qR_1(q^2)$. The pair of equations may then be solved for, in turn, $R_0(q^2)$ and $R_1(q^2)$, to give the results.
\qed
\end{proof}

To make use of the above identities, we consider bilateral sums related to the eighth order mock theta functions (the eight order equivalent of Theorem \ref{mt6t2}).
\begin{Theorem}\label{mt8t2}
(i) If $|q|<1$, then
\begin{multline}\label{28m8eq1}
G_8(a,b,d,z,q):=\sum_{n=-\infty}^{\infty}\frac{(a;q^2)_r z^r q^{r^2}}{(b;q^2)_r}
=\frac{(-zq,-qb/az;q^2)_{\infty}}{(b,q^2/a;q^2)_{\infty}}\\
\times\sum_{n=-\infty}^{\infty}\frac{(-azq/b;q^2)_r}{(-zq;q^2)_r}
\left(
-b
\right)^r q^{r^2-r}.
\end{multline}
(ii) If $|q|,  |b/a|<1$, then
\begin{equation}\label{28m8eq2}
\sum_{n=-\infty}^{\infty}\frac{(a;q^2)_r z^r q^{r^2}}{(b;q^2)_r}
=\frac{(b/a,-qb/az;q^2)_{\infty}}{(b;q^2)_{\infty}}
\times\sum_{n=-\infty}^{\infty}(-azq/b,a;q^2)_r
\left(
\frac{b}{a}
\right)^r.
\end{equation}
(iii) If $|q|, |bd/azq|<1$, then
\begin{multline}\label{28m8eq3}
\sum_{n=-\infty}^{\infty}\frac{(a;q^2)_r z^r q^{r^2}}{(b;q^2)_r}
=\frac{(-bq/az,-qz;q^2)_{\infty}}{(-q^3/az,-aqz;q^2)_{\infty}}\\
\times\sum_{n=-\infty}^{\infty}\frac{(1+azq^{4r-1})(a,-azq/b;q^2)_r\left(
baz^2
\right)^rq^{4r^2-4r}}
{(1+az/q)(b,-zq;q^2)_r}.
\end{multline}
\end{Theorem}
\begin{proof}
Let $d\to 0$ in \eqref{28m6eq1}, \eqref{28m6eq2} and \eqref{28m6eq3}, respectively.
\qed
\end{proof}

The following identities are a consequence of combining the results in Theorems \ref{mt8t2}  and \ref{st01t1}.

\begin{Corollary}\label{csto1}
If $|q|<1$, then
\begin{align}\label{cst01eq1}
\sum_{r=-\infty}^{\infty}\frac{q^{4r^2}}{(-q^2;q^2)_{2r}}
&=
\frac{(q^2;q^2)_{\infty}}{2}\left[(q;q^2)^3_{\infty}+  (-q;q^2)^3_{\infty}\right],\\
\sum_{r=-\infty}^{\infty}\frac{q^{4r^2+4r}}{(-q^2;q^2)_{2r+1}}
&=
\frac{(q^2;q^2)_{\infty}}{2q}\left[(-q;q^2)^3_{\infty}-  (q;q^2)^3_{\infty}\right],\numberthis \label{cst01eq2}\\
\sum_{r=-\infty}^{\infty}(8r+1)q^{8r^2+2r}
&=
\frac{(q^2;q^2)^3_{\infty}}{2}\left[(q;q^2)^3_{\infty}+  (-q;q^2)^3_{\infty}\right],\numberthis \label{cst01eq3}\\
\sum_{r=-\infty}^{\infty}(8r+3)q^{8r^2+6r}
&=
\frac{(q^2;q^2)^3_{\infty}}{2q}\left[(-q;q^2)^3_{\infty}-  (q;q^2)^3_{\infty}\right],\numberthis \label{cst01eq4}
\end{align}
\end{Corollary}
\begin{proof}
The identity at \eqref{cst01eq1} follows upon setting $a=-q$, $b=-q^2$ and $z=1$ in \eqref{28m8eq2}, reversing the order of summation in the resulting series on the right, replacing $q$ with $q^2$ and using \eqref{stbilatsumeq1} in conjunction with \eqref{stbilatsumprodeq1}. The identity at \eqref{cst01eq2} follows similarly, except that $z=q^2$, and \eqref{stbilatsumeq2} is used in conjunction with \eqref{stbilatsumprodeq2}.

The identities at \eqref{cst01eq3} and \eqref{cst01eq4} follow similarly from \eqref{28m8eq3}.
 For \eqref{cst01eq3}, replace $a$ with $-zq$, $b$ with $-zq^2$
 and take the limits as $z \to 1$. For \eqref{cst01eq4}, replace $z$ with $zq^2$, set $a=-zq$, $b=-zq^2$
 and again take the limits as $z \to 1$.  The details are omitted.
\qed
\end{proof}

The identities in Theorem \ref{st01t1} also contain implications for the limiting behaviour of each of the four eighth order mock theta functions that appear in these identities, as $q$ tends to certain classes of roots of unity from within the unit circle. We state these for $S_0(q)$ and $S_1(q)$, as those for $T_0(q)$ and $T_1(q)$ are equally easily derived. To avoid fractional exponents, we state the results for $q^2$ instead of $q$.

{\allowdisplaybreaks
\begin{Corollary}\label{cmt8lim}
(i) If $\zeta$ is a primitive even-order $8k$ root of unity, then, as $q$
approaches $\zeta$ radially within the unit disk, we have that
\begin{multline}\label{S0lim}
\lim_{q\to \zeta}
\left(S_0(q^2)- \frac{(q^2;q^2)_{\infty}\left[(q;q^2)^3_{\infty}+  (-q;q^2)^3_{\infty}\right]}{2(-q^2;q^2)_{\infty}}\right)\\
= -2\sum_{n=0}^{k-1} \frac{(1+\zeta^4)(1+\zeta^8)\dots (1+\zeta^{4n})}
{(1+\zeta^2)(1+\zeta^6)\dots (1+\zeta^{4n+2})}\zeta^{2n^2+6n+4}.
\end{multline}
(ii) If $\zeta$ is a primitive even-order $8k$ root of unity, then, as $q$
approaches $\zeta$ radially within the unit disk, we have that
\begin{multline}\label{S1lim}
\lim_{q\to \zeta}
\left(S_1(q^2)- \frac{(q^2;q^2)_{\infty}\left[(-q;q^2)^3_{\infty}-  (q;q^2)^3_{\infty}\right]}{2q(-q^2;q^2)_{\infty}}\right)\\
= -2\sum_{n=0}^{k-1} \frac{(1+\zeta^4)(1+\zeta^8)\dots (1+\zeta^{4n})}
{(1+\zeta^2)(1+\zeta^6)\dots (1+\zeta^{4n+2})}\zeta^{2n^2+2n}.
\end{multline}
\end{Corollary}
}

Finally, as was done for sixth order mock theta functions, we compare the results in the present paper for eighth order mock theta functions with those implied by   Mortenson's identity   at \eqref{morteq}. In that identity, if $q$ is replaced $q^8$, $x$ is set equal to $q$ and the identity of Gordon and McIntosh \cite[p. 125]{GM12},
\begin{equation}\label{mock8radeq1}
S_0(-q^2)=\frac{j(-q,q^2)j(q^6,q^{16})}{j(q^2,q^8)}-2qg_2(q,q^8)
\end{equation}
is used to replace $g_2(q,q^8)$, then
\begin{multline}\label{mock8eq2}
S_0(-q^2)- \sum _{n=0}^{\infty} \frac{q^{8 n+1}
   \left(q,q^7;q^8\right){}_n}
   {  \left(-q^8;q^8\right){}_n}
   =
\frac{q J_{1,8} }{J_{16}}
\sum_{n=0}^{\infty}
\frac{q^{8 n(n+1)}}{\left(-q,-q^7;q^8\right){}_{n+1}}\\
   -\frac{2 qJ_{16}^3 }{J_{8,16}
   J_{2,16}}
   +\frac{J_{16,32}^2
   J_{1,8}}{\bar{J}_{1,8}\bar{J}_{10,16}}
   -\frac{J_{16}^{10}
  \bar{J}_{2,16} }{J_8^4
   J_{32}^4 J_{2,16}
 \bar{J}_{10,16}}
   +\frac{\bar{J}_{1,2}
  J_{6,16} }
  {J_{2,8}}.
\end{multline}
As with sixth order mock theta functions,  when $q$ tends to a primitive root of unity of even order, the second series on the left side of  \eqref{mock8eq2} does not terminate, so that the usual kind of explicit radial limit is not obtained. However, when $\zeta$ is a primitive root of a certain order, what is obtained is a \emph{convergent} infinite series, thus leading to another type of explicit radial limit. For example, if
\[
\zeta_8=e^{2\pi i/8}=\frac{\sqrt{2}+\sqrt{2}i}{2},
\]
a primitive eighth root of unity, then it follows from \eqref{mock8eq2} that
\begin{multline}\label{mock8radeq3}
\lim_{q\to\zeta_8}
\bigg(S_0(-q^2)- \bigg[
   -\frac{2 qJ_{16}^3 }{J_{8,16}
   J_{2,16}}
   -\frac{J_{16}^{10}
  \bar{J}_{2,16} }{J_8^4
   J_{32}^4 J_{2,16}
 \bar{J}_{10,16}}
   +\frac{\bar{J}_{1,2}
  J_{6,16} }
  {J_{2,8}}\bigg ] \bigg )\\
=  \sum _{n=0}^{\infty} \frac{\zeta_8^{8 n+1}
   \left(\zeta_8,\zeta_8^7;q^8\right){}_n}
   {  \left(-\zeta_8^8;\zeta_8^8\right){}_n}
 =  \zeta_8\sum _{n=0}^{\infty}\left(\frac{(1-\zeta_8)(1-\bar{\zeta_8})}{2} \right)^n\\
   =\zeta_8\sum _{n=0}^{\infty}\left(\frac{2-\sqrt{2}}{2} \right)^n
   =1+i.
\end{multline}
Note how this compares with the radial limit given by \eqref{S0lim}:
\begin{multline}\label{S0lim2}
\lim_{q\to \zeta_8}
\left(S_0(-q^2)-\frac{\left(\left(-i q;-q^2\right)_{\infty
   }^3+\left(i q;-q^2\right)_{\infty
   }^3\right)
   \left(-q^2;-q^2\right){}_{\infty }}{2
   \left(q^2;-q^2\right){}_{\infty }}\right)\\
= 1+i.
\end{multline}
While the limits are the same, the two theta functions subtracted from $S_0(-q^2)$ are not equal as functions of $q$.

Here also, as with the sixth order mock theta functions $\psi(q)$ and $\lambda(q)$, a radial limit is not obtained $q$ tends to a primitive root of unity of odd order. What is true is that if $\zeta$ is a primitive root of unity of order $2k+1$, then
\begin{multline}\label{mock8radeq4}
\sum _{n=0}^{k} \frac{(-1)^n\zeta^{2n^2}
   \left(\zeta^2;\zeta^4\right){}_n}
   {  \left(-\zeta^4;\zeta^4\right){}_n}- \sum _{n=0}^{k} \frac{\zeta^{8 n+1}
   \left(\zeta,\zeta^7;\zeta^8\right){}_n}
   {  \left(-\zeta^8;\zeta^8\right){}_n}\\
   =\lim_{q\to \zeta}\left(
 \frac{\bar{J}_{1,2}
  J_{6,16} }
  {J_{2,8}}
    -\frac{2qJ_{16}^3 }{J_{8,16}
   J_{2,16}}
   -\frac{J_{16}^{10}
  \bar{J}_{2,16} }{J_8^4
   J_{32}^4 J_{2,16}
 \bar{J}_{10,16}}
   \right),
\end{multline}
and once again, experiment seems to suggest that each side is identically zero when $\zeta$ is any primitive root of unity of odd order. As with the hypotheses suggested by experiment for mock theta functions of sixth order, we have not attempted to prove these assertions.

Results similar to those described above may be derived for other eighth order mock theta functions.

\setcounter{equation}{0}

 \allowdisplaybreaks{

}
\end{document}